\documentclass[12pt]{article}

\usepackage{amsmath}
\usepackage{amsfonts}
\usepackage{amssymb}
\usepackage{graphicx, color}
\usepackage{hyperref}
\usepackage{epsfig}
\usepackage{epstopdf} 

\definecolor{darkgreen}{rgb}{0,0.55,0}

\newtheorem{proposition}{Proposition}[section]
\newtheorem{theorem}{Theorem}[section]
\newtheorem{lemma}[theorem]{Lemma}

\newtheorem{remark}[theorem]{Remark}

\DeclareSymbolFont{AMSb}{U}{msb}{m}{n}
\DeclareMathSymbol{\N}{\mathbin}{AMSb}{"4E}
\DeclareMathSymbol{\Z}{\mathbin}{AMSb}{"5A}
\DeclareMathSymbol{\R}{\mathbin}{AMSb}{"52}
\DeclareMathSymbol{\Q}{\mathbin}{AMSb}{"51}
\DeclareMathSymbol{\I}{\mathbin}{AMSb}{"49}

\newcommand{\calC}{{\mathcal C}}
\newcommand{\calS}{{\mathcal S}}

\numberwithin{equation}{section}

\begin{document}

\title{ The Sphere Covering Inequality and Its Applications}

\author{{Changfeng Gui\footnote{Department of Mathematics, University of Texas at San Antonio, Texas, USA. E-mail: changfeng.gui@utsa.edu.}
\qquad Amir Moradifam\footnote{Department of Mathematics, University of California, Riverside, California, USA. E-mail: moradifam@math.ucr.edu. }}}
\date{\today}

\smallbreak \maketitle

\begin{abstract}
In this paper, we show that   the  total area  of two {\it distinct}  surfaces with Gaussian curvature  equal to 1,   which  are  also conformal to  the Euclidean unit disk  with the same conformal factor on the boundary,  must be at least $4 \pi$.  In other words,  the areas of these surfaces must cover the whole unit sphere after a proper rearrangement.   We refer  to this lower bound of total area as the Sphere Covering Inequality.  The inequality  and its generalizations  are applied to a number of open problems related to  Moser-Trudinger  type inequalities,  mean field equations and Onsager vortices, etc,  and yield   optimal results. 

\end{abstract}
\maketitle

\section{Introduction} 

A large number of  important second order nonlinear  elliptic equations  involve  exponential nonlinearities.  These equations arise, for example,  in 
the study of Gaussian curvature of surfaces with metrics conformal to Euclidean metric ( \cite{CLiu-MR1209959}, \cite{CY2-MR925123}, \cite{CY1-MR908146},  \cite{CL1-MR1121147},\cite{CL2-MR1338474}, etc.), Moser-Trudinger type inequalities (\cite{Aubin-MR0448404}, \cite{Aubin2-MR534672}, \cite{Beckner-MR1230930}, \cite{DEJ-MR3377875}),  \cite{DET-MR2509375}, \cite{GL-MR2670931}, \cite{GW-MR1760786}, \cite{Moser-MR0301504}, \cite{O-MR677001},  \cite{OPS-MR960228}, \cite{T-MR0216286}), the mean field theory of statistical mechanics of  classical vortices  and thermodynamics (\cite{ BD-MR3355004}, \cite{BLT-MR2838340}, \cite{CLMP1-MR1145596},  \cite{CLMP2-MR1362165},  \cite{CCL-MR2055839}, \cite{CK-MR1262195}, \cite{K-MR1193342}, \cite{Lin1-MR1770683}),and  self gravitating cosmic string 
configurations in the framework of Einstein's general relativity  (\cite{CGS-MR2944090}, \cite{PT-MR2876669}, \cite{Yang-MR1277473}).
In this article,  we shall prove a  basic and important inequality which  becomes a crucial tool for tackling several open problems in  the above
mentioned areas. 

Let us consider the equation
\begin{equation}\label{model}
\Delta v +e^{2v}=0, \quad y \in \Omega,
\end{equation}
where $\Omega \subset \R^2$ is a $C^2$ simply-connected bounded region.  It is well-known that for a solution $v \in C^2(\bar \Omega)$ of \eqref{model},  the  two dimensional Riemannian  manifold with boundary $(\Omega, g)$  with a conformal Euclidean metric $dg= e^{2v} dy$   has Gaussian curvature equal to 1 everywhere.   The total area as well as the total curvature of such manifold is equal to  $A =\int_{\Omega} e^{2v}dy$.  The well-known Gauss-Bonnet Theorem states that
$$
A=\int_{\Omega} e^{2v}dy=  \int_{\Omega} dg=  2\pi -\int_{\partial \Omega} \kappa_g dg_s
$$
where $\kappa_g$ is the geodesic curvature  and $dg_s$  the length parameter of $\partial \Omega$.  
From the equation,  it is also easy to see that
$$
A=-\int_{\partial \Omega} \frac{\partial v}{\partial r} ds. 
$$
These formulas, though very useful in general,  do not impose any restriction  on the area of the surface,  as  the uniformization theorem says that every simply-connected  Riemann surface is conformally equivalent to one of the three domains: the open unit disk, the complex plane, or the Riemann sphere.   However,  if there is another surface $(\Omega, \tilde g)$ with a distinct conformal metric $d \tilde g= e^{2\tilde v} dy$ in $\Omega$,  where  $\tilde v \in C^2(\bar \Omega)$  is a solution of  \eqref{model} and $\tilde g=g $ on $\partial \Omega$,  we shall show 

\begin{equation}\label{key}
\tilde A + A = \int_{\Omega} (e^{2 \tilde v} + e^{2v}) dy \ge 4 \pi.
\end{equation}

Since the standard sphere has Gaussian curvature 1 and area $4 \pi$, and  these two surfaces have total area bigger than  or equal to that of  the standard sphere,  one may think that  these two surfaces could cover the standard sphere if  they are properly arranged (this will be made more rigorous later in Section 2.1).   The equality obviously hold when the two surfaces are isometric to two complementing spherical caps
on the standard sphere.  We therefore refer to \eqref{key}  as the Sphere Covering Inequality.  Let us give an explicit example to better demonstrate the Sphere Covering Inequality. Set $\lambda>0$ and define 
\begin{equation}\label{VLambda}
V_{\lambda}(y):=-\ln(1+\frac{\lambda^2 |y|^2 }{8})+\ln(\lambda)-\frac{1}{2}\ln(2),
\end{equation}
satisfying
\begin{equation*}
\Delta V_\lambda+e^{2 V_\lambda}=0, \quad y \in \R^2.
\end{equation*}
For every $\lambda_2>\lambda_1$, there exists a unique $R \in \R$ such that $V_{\lambda_1}(\partial B_R)=V_{\lambda_2}(\partial B_R)$ and $V_{\lambda_2}> V_{\lambda_1}$ in $B_R$.  
Now define two surfaces ${\cal S}_1$ and ${\cal S}_2$ with constant Gaussian curvature 1 as follows
\[ {\cal S}_{1}=\bigl(B_R, e^{2V_{\lambda_1}}dy\bigr) \ \ \hbox{ and}
 \ \  {\cal S}_{2}=\bigl(B_R, e^{2V_{\lambda_2}}dy\bigr).\]
Notice that the metrics $g_{i}=e^{2V_{\lambda_i}}dy, i=1, 2$ have the same conformal factor on $\partial B_R$, which 
implies that $\lambda_1 \lambda_2=\frac{8}{R^2}$.  By scaling in $y \in \R^2$ and  using the stereographic projection $\Pi$: $S^2\rightarrow \R^2$ with respect to the north pole $N=(1,0,0)$: 
\[y=\Pi(x):= \left( \frac{x_1}{1-x_3}, \frac{x_2}{1-x_3}\right),\]
we can see that the surface  
${\cal S}_{1}=(B_R, e^{2V_{\lambda_1}}dy)$ is isometric to 
$$
\bigl(B_{\lambda_1R/\sqrt{8}}, e^{2V_1}dy\bigr)=\bigl(B_{\lambda_1R/\sqrt{8}}, \frac{1}{(1+|y|^2)^2} dy\bigr),
$$
which is isometric to a disc $\calC_1$ around the south pole.  Similarly  the surface 
${\cal S}_{2}=(B_R, e^{2V_{\lambda_1}}dy)$ is isometric to  
$$
\bigl(B_{\lambda_2R/\sqrt{8}}, e^{2V_1}dy\bigr)=\bigl(B_{\lambda_2R/\sqrt{8}}, \frac{1}{(1+|y|^2)^2} dy\bigr),
$$
which is isometric to  a disc $\calC_2$ around the south pole.

Using the Kelvin transformation $z= y/|y|^2$ and the fact $\lambda_1 \lambda_2=\frac{8}{R^2}$,  one can see that 
$$
\bigl(B_{\lambda_2R/\sqrt{8}}, \frac{1}{(1+|y|^2)^2} dy),
$$
is isometric to 
$$
\bigl(\R^2 \setminus B_{\lambda_1R/\sqrt{8}}, \frac{1}{(1+|z|^2)^2} dz \bigr),
$$
which is isometric to $ S^2 \setminus \calC_1$. 
This implies  that  ${\cal S}_{1}$ and ${\cal S}_{2}$ are indeed  isometric to two complimenting spherical caps on the unit sphere, and therefore their total areas are exactly $4 \pi.$

We will prove the inequality (\ref{key}) in a more general setting as follows.

\begin{theorem}[The Sphere Covering Inequality]\label{MainTheorem0}
Let $\Omega $ be a simply-connected  subset of $\R^2$ and assume $v_i \in C^2(\overline{\Omega})$, $i=1,2$ satisfy
\begin{equation}\label{mainpde}
\Delta v_i +e^{2v_i}=f_{i}(y),
\end{equation}
where $f_2 \geq f_1 \ge 0$  in $\Omega$. If $v_2 \ge v_1,  v_2 \not \equiv v_1$ in $\Omega$ and $v_2=v_1$ on $\partial \Omega$, then
\begin{equation}
\int_{\Omega} e^{2v_1}+e^{2v_2}dy \geq 4\pi. 
\end{equation}
Moreover,  the equality only holds when $f_2 \equiv f_1 \equiv 0$ and $(\Omega, e^{2v_i}dy ) $, $i=1,2$ are isometric to  two complimenting spherical caps on the standard unit sphere.  
\end{theorem}

For the simplicity of the equation,  we may replace $2v$ by $u-\ln2$ and consider

\begin{equation}\nonumber
\Delta u +e^{u}=0, \quad y \in \Omega,
\end{equation}

Geometrically, this is equivalent to  multiplying the conformal factor by  $\sqrt 2$ so the sphere in comparison 
has radius $\sqrt 2$ and total area $8 \pi$. Indeed Theorem \ref{MainTheorem0} is equivalent to Theorem \ref{MainTheorem} in Section 3.

The Sphere Covering Inequality is closely related to the symmetry of solutions of elliptic equations with exponential nonlinearity in $\R^2$. To see the connection, consider the equation 
\begin{equation}\label{SCIpde}
\Delta w+e^w=f \geq 0 \ \  \ \ \hbox{in}\ \ \R^2,
\end{equation}
and let $w$ be a classical solution with a critical point located at $P\in \R^2$.  Assume that $f$ is evenly symmetric about a line passing through
$P$.   It follows from the Sphere Covering Inequality that if $\int_{\R^2}e^w dy< 16 \pi$, then $u$ must be symmetric about the line.  More precisely, suppose $P=(p,0)$  and $f(y_1, y_2)=f(y_1, -y_2)$ in $\R^2$.  Define $\bar{w}(y_1,y_2)=w(y_1,-y_2)$, and set 
\[\tilde v:=w-\bar{w}.\]
Then $\tilde v$ satisfies 
\[\Delta \tilde v+c \tilde v=0,\]
where 
\[c=\frac{e^{w}-e^{\bar{w}}}{w-\bar{w}}.\]
Suppose $\tilde v\not \equiv 0$.  It follows from the Hopf's lemma that the nodal line of $\tilde v$ divides a neighborhood of $P$ into at least four regions, and consequently there exist at least two simply-connected regions $\Omega_1, \Omega_2 \subset \R^2_{+}$ such that $\tilde v>0$ in $\Omega_1$, $\tilde v<0$ in $\Omega_2$, and $\tilde v=0$ on $\partial \Omega_1 \cup \partial \Omega_2$.  Therefore on each $\Omega_i$, $i=1,2$, the equation \eqref{SCIpde} has two distinct solutions, $w$ and $\bar{w}$, satisfying the assumptions of Theorem \ref{MainTheorem}, which is  an equivalent form of Theorem \ref{MainTheorem0} for \eqref{SCIpde}. Thus 
\[\int_{\R^2}e^{w}dy \geq \int_{\Omega_1}e^{w}+e^{\bar{w}}dy+\int_{\Omega_2}e^{w}+e^{\bar{w}}dy\geq 16\pi,\]
which is a contradiction and leads to the symmetry of $w$. 

The above argument is at the core of the proof of the symmetry results in this paper,  and  consists of two main ingredients:  the Sphere Covering
Inequality and the nodal set analysis.  The idea of using nodal sets to prove symmetry results for elliptic equations with exponential nonlinearity  was used by Lin and others to obtain symmetry results for mean field equations in $\R^2$ and on $S^2$ and flat tori (see, e.g.,  \cite{Lin2-MR1781481},  \cite{LL-MR2352519},  \cite{GL-MR2670931}, \cite{BLT-MR2838340}, etc).   The key in their arguments is  Proposition \ref{4piBound} in Section 3,  which has a  geometrical interpretation in terms of the extremal first eigenvalue (see Remark \ref{Remark}) and usually yields  a lower bound  of $2 \pi$ for \eqref{mainpde} instead of $4 \pi$ in the Sphere Covering Inequality.  Note that Proposition  \ref{4piBound} may be regarded as a special limiting case of the Sphere Covering Inequality,  although  the  geometric meaning of the Sphere Covering Inequality itself  still remains unclear and worth further exploring.  In this paper,  the Sphere Covering Inequality will be used to solve several important open problems. Below we introduce some of the problems.

\subsection{Best Constant in A Moser-Trudinger Type  Inequality}
Let $S^2$ be the unit sphere and for $u\in H^1(S^2)$ define 

\begin{equation}\label{MTAO}
J_\alpha(u)=\frac{\alpha}{4} \int_{S^2}|\nabla u| ^2 d\omega+ \int_{S^2}u d\omega-\log \int_{S^2}e^{u}d \omega,
\end{equation}
where the volume form $d\omega$ is normalized so that $\int_{S^2}d\omega=1$. The well-known Moser-Trudinger inequality \cite{Moser-MR0301504} says that $J_\alpha$ is bounded below if and only if $\alpha \geq 1$. Onofri \cite{O-MR677001} showed that for $\alpha \geq 1$ the best lower bound is equal to zero. On the other hand, Aubin \cite{Aubin2-MR534672} proved that if $J_\alpha$ is restricted to  
\[\mathcal{M}:=\{u\in H^1(S^2): \ \ \int_{S^2}e^ux_i=0, \ \ i=1,2,3\},\]
then for $\alpha \geq\frac{1}{2}$, $J_\alpha$ is bounded below and the infimum is attained in $\mathcal{M}$. In 1987 Chang and Yang \cite{CY2-MR925123}, in their work on prescribing Gaussian curvature on $S^2$ (see also \cite{CY1-MR908146}), showed that for $\alpha$ close to 1 the best constant again is equal to zero, and this led to the following conjecture. \\ \\
{\bf Conjecture A.} For $\alpha \geq \frac{1}{2}$ 
\[\inf_{u\in \mathcal{M}} J_{\alpha}(u)=0.\]

In 1998,  Feldman, Froese, Ghoussoub and the first author \cite{FFGG-MR1606461} proved that this conjecture is true for axially symmetric functions when $\alpha>\frac{16}{25}-\epsilon$. Later the first author and Wei \cite{GW-MR1760786}, and independently Lin \cite{Lin1-MR1770683} proved Conjecture A for radially symmetric function, but the problem remained open for non-axially symmetric functions. 

In \cite{GL-MR2670931}  Ghoussoub and Lin, showed that Conjecture A holds true for $\alpha \geq \frac{2}{3}-\epsilon$, for some $\epsilon>0$. See Chapter 19 in  \cite{GoussoubMoradifam} for a complete history of the problem. In this paper, among other results, we will prove that  Conjecture A is true. 

\begin{theorem}\label{MAOTheorem}
For $\alpha \geq \frac{1}{2}$ 
\[\inf_{u\in \mathcal{M}} J_{\alpha}(u)=0.\]
\end{theorem} 
\vspace{.8cm}
Indeed we apply Theorem \ref{MainTheorem} to show that the corresponding Euler-Lagrange equation 
\begin{equation}
\frac{\alpha}{2}\Delta u+\frac{e^u}{\int_{S^2}e^u d\omega}-1=0 \ \ \hbox{on}\ \ S^2
\end{equation}
has only constant solutions. 

\subsection{A Mean Field Equation with singularity on $S^2$}
Consider the mean field equation 
\begin{equation}\label{meanFieldPDE}
\Delta_g u+\lambda \left( \frac{e^{u}}{\int_{S^2}e^u d\omega}-\frac{1}{4\pi}\right)=4\pi(\delta(P)-\frac{1}{4\pi}) \ \ \hbox{on} \ \ S^2, 
\end{equation}
where $g$ is the standard metric on $S^2$ with the corresponding volume form $d\omega$, $\alpha>-1$, $\lambda >0$, and $P\in S^2$. In \cite{BLT-MR2838340}, motivated by the study of  vortices in self-dual gauge field theory,  Bartolucci, Lin and Tarentello studied  symmetry of solutions of (\ref{meanFieldPDE}) under the assumption 
\begin{equation}\label{meanFieldPDE1}
\lambda=4\pi(3+\alpha)
\end{equation}
and showed that (\ref{meanFieldPDE}) admits a solution if and only if $\alpha \in (-1,1)$. Then, via a new bubbling phenomenon, they proved that there exists $\delta>0$ such that for $\alpha \in (1-\delta,1)$ the equation (\ref{meanFieldPDE})  admits a unique solution that is in addition is axially symmetric about the direction $\overrightarrow{OP}$, and proposed the following  \\ \\
{\bf Conjecture B.}  All solutions of (\ref{meanFieldPDE})-(\ref{meanFieldPDE1}) are axially symmetric about $\overrightarrow{OP}$ for every $\alpha \in (-1,1)$.\\ \\
In Section 5,  we shall use the Sphere Covering Inequality to provide an affirmative answer to the above question. Indeed we will prove the following result.

\begin{theorem}\label{MeanFieldUniquenessTheorem} For every $\alpha \in (-1,1)$  equation (\ref{meanFieldPDE})-(\ref{meanFieldPDE1}) has a unique solutions that in addition is axially symmetric about $\overrightarrow{OP}$.  
\end{theorem}

\subsection{A Mean Field Equation for the Spherical Onsager Vortex}

Consider the following equation 
\begin{equation}\label{OnsagerVortexPDE}
\Delta_g u(x)+\frac{\exp(\alpha u(x)-\gamma \langle n, x \rangle )}{\int_{S^2} \exp(\alpha u(x)-\gamma \langle n, x \rangle) d\omega}-\frac{1}{4\pi}=0 \ \ \hbox{on}\ \ S^2,
\end{equation}
where $g$ is the standard metric on $S^2$ with the corresponding volume form $d\omega$, $ \vec{n}$ is a unique vector in $R^3$, $\alpha\geq 0$, and $\gamma \in R$. Since $\gamma <0$ can be changed to $-\gamma$ by replacing the north pole with the south pole,  we only need to consider the case $\gamma \ge 0$.  This equation is invariant up to adding a constant and we seek a normalized solution
with 
\begin{equation}
\int_{S^2}u d\omega=0.
\end{equation}
In \cite{ Lin1-MR1770683},  Lin showed that if $\alpha< 8\pi$, then for $\gamma \geq 0$ the equation (\ref{OnsagerVortexPDE}) has a unique solution  that in addition is axially symmetric with respect to $\vec{n}$.  In this case  the coefficient in the equation is  decreasing and therefore the moving plane method applies (see (\ref{onsager})).  He also conjectured the following \\ \\
{\bf Conjecture C. } Let $\gamma\geq 0$ and  $\alpha \leq 16\pi$. Then every solution $u$ of (\ref{OnsagerVortexPDE}) is axially symmetric with respect to $\vec{n}$. \\ 

In an attempt to prove this conjecture, in \cite{Lin2-MR1781481}, C.S. Lin proved the following theorems for $\alpha > 8\pi$. \\ \\
{\bf Theorem A.} (\cite{Lin2-MR1781481}) For every $\gamma >0$, there exists $\alpha_0=\alpha_0(\gamma)>8 \pi$ such that, for $8\pi < \alpha \leq \alpha_0$, any solution $u$ of (\ref{OnsagerVortexPDE}) is axially symmetric. \\ \\
{\bf Theorem B.} (\cite{Lin2-MR1781481}) Let $u_i$ be a solution of (\ref{OnsagerVortexPDE}) with $\gamma=0$ and $\alpha_i \rightarrow 16 \pi$. Suppose $\lim_{i \rightarrow \infty} \sup u_i(x)=+\infty$. Then $u_i$ is axially symmetric with respect to some direction $ \vec{n}_i$ in $R^3$  for $i$ large enough. \\

In Section 5,  we shall apply the Sphere Covering Theorem to prove the following result. 
\begin{theorem}\label{OnsagerVortexTheorem}
Suppose $8 \pi \leq \alpha \leq 16 \pi$ and 
\begin{equation}
0 \leq \gamma \leq \frac{\alpha}{8\pi}-1. 
\end{equation}
Then every solution of (\ref{OnsagerVortexPDE}) is axially symmetric with respect to $\vec{n}$. In particular if $\gamma=0$ and $8 \pi \leq \alpha \leq 16 \pi$, then the trivial solution $u\equiv 0$ is the only solution of (\ref{OnsagerVortexPDE}).
\end{theorem}

In all  problems above and many others,  there exists a critical number  $8 \pi$ for a quantity which may be interpreted as 
total area literally.  In \eqref{MTAO},   the quantity is  $ 4 \pi /\alpha$,  with  $ \alpha$ being  a parameter;  In \eqref{meanFieldPDE},  the quantity is 
$\lambda$; While in \eqref{OnsagerVortexPDE},  the quantity is $\alpha$.  Note that  the parameter $\alpha $ has different meanings in these three equations  which should be clear from the context.  The work of Brezis and Merle \cite{BM-MR1132783} and Y.Y. Li \cite{YLi-MR1673972} as well as others showed that $8 \pi m,  m \in N$ are values where  solutions of these type of equations may lose compactness and blow-up 
phenomena may happen.  The critical level $8 \pi$ also separates two  significantly different cases in terms of the coerciveness of associate functionals and the positiveness of linearized operators. A crucial tool is required to deal with the supercritical cases of many important problems in 
related research.   The Sphere Covering Inequality provides exactly such a  much needed tool.   Besides the applications in this paper,   other
applications  of the Sphere Covering Inequality will also be discussed in  forthcoming papers \cite{GM-Torus} and \cite{GM-MeanField}. 

  The paper is organized as follows.   In Section 2,  we shall discuss some preliminary results about the classical Bol's inequality and prove a counterpart of Bol's inequality for radially symmetric functions which is needed for the proof of the Sphere Covering Inequality.  In Section 3, we will prove  the
  Sphere Covering Inequality.   In Section 4,   the Sphere Covering Inequality shall be applied to \eqref{MTAO} to show the best constant. 
  Finally,  we will present a general symmetry result regarding Gaussian curvature equations on $\R^2$ which leads to  optimal results for
  \eqref{meanFieldPDE}, \eqref{OnsagerVortexPDE} and others.

\section{Bol's Isoperimetric Inequality}
Bol's isoperimetric inequality plays a crucial role in the proof of our main results. In this section we present some preliminary results on Bol's inequality that will be used in subsequent sections. First we recall  the classical  Bol's isoperimetric inequality \cite{Bandle-MR572958, Bol-MR0018858, Suzuki-MR1186683}: 

\begin{proposition}\label{MainBolProp}
Let $\Omega \subset \R^2$ be a simply-connected and assume $u\in C^2(\Omega)$  satisfies 
\begin{equation}\label{MainPDECond}
\Delta u+e^{u}\geq 0 \ \ \hbox{and} \ \ \int_{\Omega}e^{u}\leq 8 \pi.
\end{equation}
Then for every $\omega \Subset \Omega$ of class $C^1$ the following inequality holds
\begin{equation}\label{BolInequality}
\left( \int_{\partial \omega}e^{\frac{u}{2}}\right)^2\geq \frac{1}{2}\left( \int_{\omega} e^u\right) \left( 8 \pi-\int_{\omega}e^u\right).
\end{equation}
\end{proposition}

We first show an example for which the equality in Bol's inequality holds. 
For $\lambda>0$,   let define $U_\lambda$   by 
\begin{equation}\label{ULambda}
U_{\lambda}:=-2\ln(1+\frac{\lambda^2 |y|^2 }{8})+2\ln(\lambda).
\end{equation}
Then
\begin{equation*}
\Delta U_\lambda+e^{U_\lambda}=0,
\end{equation*}
and 
\[\int_{B_{r}}e^{U_\lambda}dy= \frac{8 \pi \lambda^2r^2}{8+\lambda^2r^2},\]
for all $r>0$, where $B_r$ denotes the ball of radius $r$ centered at the origin in $\R^2$. 
One can check that 
\[\left( \int_{\partial B_r}e^{\frac{ U_\lambda}{2}}\right)^2= \frac{1}{2}\left( \int_{B_r} e^{ U_\lambda}\right) \left( 8 \pi-\int_{B_r}e^{ U_\lambda}\right),\]
for all $r>0$ and $\lambda>0$.  Indeed, $e^ {U_{\lambda} (y)}dy$ corresponds to the metric in a standard sphere with radius $\sqrt 2$.

By examining the proof of Bol's inequality (see, e.g., \cite{Suzuki-MR1186683}), it can be seen that if the equality holds for some  $\omega$
in \eqref{BolInequality}, then  $\Delta u+e^{u} = 0 $ in $\omega$,  and $ e^{u(z) }dz= e^{U_\lambda(\xi)}d\xi $, where $z=g^{-1}(\xi)$ for some analytic function $g: \Omega \rightarrow B_R$, and $\lambda>0$. More precisely, let us consider the case that $\omega=\Omega$ is simply-connected , and  follow the arguments in \cite{Suzuki-MR1186683}
by considering  the  harmonic lifting  $h$ of  boundary value of $ u$  in $\Omega$, i.e.,
$$
\Delta h(z)=0, \quad z \in \Omega;  \quad h=u, \quad z \in \partial \Omega.
$$
It is known that there is an analytic function $\xi=g(z)$ such that $e^h= |g'(z)|^2$.  The equality  in \eqref{BolInequality} implies  that 
$g(\Omega)=B_R$ and $g(\partial \Omega) = g(\partial B_R)$ for some $R>0$.  Furthermore,  letting 
$$
 q(\xi): = e^{u \bigl(g^{-1}(\xi) \bigr) }  |g'  \bigl( g^{-1}(\xi) \bigr) |^{-2},
$$
we have 
$ e^u(z)dz=  q(\xi) d\xi $ and  $ q(\xi)$ is radially symmetric.   Therefore $q(\xi)=e^ {U_{\lambda} (\xi)} $ for some $\lambda>0$.   In general,  if 
$\omega $ is not simply-connected ,  Bartolucci and Lin  showed in \cite{BL-MR3201892} that  strict inequality holds in  \eqref{BolInequality}.
So  the equality may hold only for simply-connected   $\omega$.  
Note that if $\Omega$ is not simply-connected,  Proposition \ref{MainBolProp} is not valid in general.  Indeed,   \eqref{BolInequality}  does not hold for certain annulus regions  as shown in \cite{BL-MR3201892}.

For the proof of our main results, we shall need the following counterpart of the Bol's inequality for radial functions. The proof is a modification of an argument by Suzuki \cite{Suzuki-MR1186683} and we present it here for the sake of completeness. 

\begin{proposition}\label{RadialBolProp}
Let $B_R$ be the ball of radius $R$ in $\R^2$ $\psi \in C^{0,1} (\overline{B_R})$ be a strictly decreasing, radial,  Lipschitz function satisfying
\begin{equation}\label{MainPDECondRelaxed}
\int_{\partial B_r}|\nabla \psi |ds \le \int_{B_r}e^{\psi}  dy \ \ \hbox{a.e. }  r\in(0,R), \ \ \hbox{and} \ \ \int_{B_R}e^{\psi}\leq 8 \pi.
\end{equation}
Then the following inequality holds
\begin{equation}\label{RadialBolInequality}
\left( \int_{\partial B_R}e^{\frac{\psi}{2}}\right)^2\geq  \frac{1}{2}\left( \int_{B_R} e^\psi\right) \left( 8 \pi-\int_{B_R}e^\psi \right).
\end{equation}
Moreover if $\int_{\partial B_r}|\nabla \psi |ds  \not \equiv \int_{B_r}e^{\psi }  dy $  in $(0,R)$, then the inequality in (\ref{RadialBolInequality}) is strict. 
\end{proposition}

{\bf Proof.} Let $\beta:=\psi(R)$ and define
\[k(t)=\int_{\{\psi>t\} }e^{\psi}dy, \ \ \hbox{and} \ \ \mu(t)=\int_{\{\psi>t\} }dy,\]
for $t>\beta$. Then 
\[-k'(t)=\int_{\{\psi=t\}}\frac{e^{\psi}}{|\nabla \psi|}=-e^{t}\mu'(t), \ \ \hbox{for  \it {a.e.}}  \ \ t>\beta.\]
Hence 
\begin{eqnarray}\label{RelaxedBolSharp}
-k(t)k'(t)&\geq & \int_{\{\psi=t\} }|\nabla \psi| \cdot \int_{\{\psi=t\}}\frac{e^{\psi}}{|\nabla \psi|}\\
&= & (\int_{\{\psi=t\}}e^{\psi/2})^{2}=e^{t}(\int_{\{\psi=t\}}ds)^{2}\nonumber \\
&= &4 \pi e^t\int_{\{\psi>t\}}dy=4 \pi e^t\mu(t), \ \ \hbox{for  \it {a.e.}}  \ \ t>\beta.  \nonumber
\end{eqnarray}
Therefore 
\[\frac{d}{dt} [e^t \mu(t)-k(t)+\frac{1}{8\pi}k^2(t)]=e^t \mu(t)+\frac{1}{4\pi}k'(t)k(t)\le 0, \ \ \hbox{for  \it {a.e.}}  \ \ t>\beta.\]
 Integrating on $(\beta, \infty)$ we get 
\begin{equation}\label{calInequality}
\left[ e^t \mu(t)-k(t)+\frac{1}{8\pi}k^2(t)\right]_{\beta}^{\infty}=-\left(e^{\beta} \mu(\beta)-k(\beta)+\frac{1}{8\pi}k^2(\beta) \right) \le 0.
\end{equation}
Now notice that 
\[ k(\beta)=\int_{B_R}e^{\psi}dy\]
and
\[e^{\beta}\mu(\beta)=e^{\beta} \int_{ B_R}dy= \frac{1}{4\pi} e^{\beta}(\int_{\partial B_R}ds)^{2}=\frac{1}{4 \pi}(\int_{\partial B_R} e^{\frac{\psi}{2}}ds)^{2}. \]
Thus (\ref{RadialBolInequality}) follows from the inequality (\ref{calInequality}). 
Finally if $\int_{\partial B_r}|\nabla \psi |ds  \not \equiv  \int_{B_r}e^{\psi}  \in(0,R)$, then the inequality (\ref{RelaxedBolSharp}) will be strict
in a set  with a positive measure in $ \{t>\beta\}$, and consequently   \eqref{calInequality}  and (\ref{RadialBolInequality}) will also be strict. \hfill $\Box$\\ 

\subsection{Rearrangement with respect to two measures}
Let $\Omega \subset \R^2$ and $\lambda>0$, and suppose that $w\in C^2(\overline{\Omega})$ satisfies
\[\Delta w+e^{w}\geq 0.\]
Then any function $ \phi \in C^2(\overline{\Omega})$ can be equimeasurably rearranged with respect to the measures $e^wdy$ and $e^{U_\lambda}dy$ (see \cite{Bandle-MR572958}, \cite{Suzuki-MR1186683}, and \cite{LL-MR2352519}), where $U_\lambda$ is defined in (\ref{ULambda}). More precisely, for $t>\min_{y\in \overline\Omega}\phi$ define
\[\Omega_t:=\{\phi>t\} \subset\subset \Omega,\]
and define $\Omega^*_t$ be the ball centered at the  origin in $\R^2$ such that 
\[\int_{\Omega^*_t}e^{U_{\lambda}}dy=\int_{\Omega_t}e^{w}dy:=a(t).\]
Then $a(t) $ is  a right-continuous function,  and $\phi^*: \Omega^* \rightarrow \R$ defined by $\phi^*(y):=\sup \{t\in \R: y\in \Omega^*_t\}$ provides an equimeasurable rearrangement of $\phi$ with respect to the measure $e^{w}dy$ and $e^{U_{\lambda}}dy$, i.e. 
\begin{equation}\label{rearrang}
\int_{\{\phi^*>t\}}e^{U_{\lambda}}dy=\int_{\{\phi>t\}}e^{w}dy, \ \ \forall  t>\min_{y\in \overline \Omega}\phi.
\end{equation}

Let
\[
j(t):=\int_{\{\phi >t \} } |\nabla \phi |^2 dy, \quad j^*(t):=\int_{\{\phi >t \} } |\nabla \phi^* |^2 dy,  \ \ \forall  t>\min_{y\in \overline \Omega}\phi;
\]
\[J(t):=\int_{\{\phi>t\}}|\nabla \phi| dy, \quad J^*(t):=\int_{\{\phi^*>t\}}|\nabla \phi^*| dy, \ \ \forall  t>\min_{y\in \overline \Omega}\phi.
\]

It is easy to see that both $j(t)$ and $J(t)$ are  absolutely continuous  and decreasing  in $ t>\min_{y\in \overline \Omega}\phi $.  

When $ \phi \equiv C $ on $\partial \Omega$,  it
can  be shown that 
\begin{equation}
\int_{\{\phi=t\}}|\nabla \phi| ds\geq \int_{\{\phi^*=t\}}|\nabla \phi^*|ds, \quad  \hbox{for  \it {a.e.}}  \ \  t>\min_{y\in \overline \Omega}\phi.
\end{equation}
Indeed it follows from  Cauchy-Schwarz and Bol's inequalities that 

 \begin{eqnarray*}
 \int_{\{\phi=t\}}|\nabla \phi| ds
&\geq & \left(\int_{\{\phi=t\}} e^{\frac{w}{2}} \right)^{2} \left( \int_{\{\phi=t\}} \frac{e^{w}}{|\nabla \phi|}\right)^{-1}\\
&=& \left(\int_{\{\phi=t\}} e^{\frac{w}{2}} \right)^{2} \left( -\frac{d}{dt}\int_{\Omega_t}e^{w}\right)^{-1}\\
&\geq &\frac{1}{2}(\int_{\Omega_t}e^{w})(8\pi - \int_{\Omega_t}e^{w})(-\frac{d}{dt}\int_{\Omega_t}e^{w})^{-1}\\
&=&  \frac{1}{2}(\int_{\Omega^*_t}e^{U_{\lambda}})(8\pi - \int_{\Omega^*_t}e^{U_{\lambda}})(-\frac{d}{dt}\int_{\Omega^*_t}e^{U_{\lambda}})^{-1}\\ 
 &=& \int_{\{\phi^*=t\}}|\nabla \phi^*| ds, \ \ \hbox{for  \it {a.e.}}  \ \ t>\min_{y\in \overline \Omega}\phi.
 \end{eqnarray*}
 
 It also follows that   $j^*(t),  J^*(t)$ are  absolutely  continuous and decreasing  in $t>\min_{y\in \overline \Omega}\phi$,  since both  functions are right-continuous by definition and 
 \[
 0\le j^*(t)-j^*(t-0)\le j(t) -j(t-0)=\int_{\{\phi=t\}}|\nabla \phi|^2 dy =0,   \ \ t>\min_{y\in \overline \Omega}\phi. \]
 \[
 0\le J^*(t)-J^*(t-0)\le J(t) -J(t-0)=\int_{\{\phi=t\}}|\nabla \phi|dy =0,   \ \ t>\min_{y\in \overline \Omega}\phi. \]

 Therefore we have  the following proposition. 
 
\begin{proposition}\label{rearrangProp}
Let $w\in C^2(\overline{\Omega})$ satisfy 
\[\Delta w+e^{w}\geq 0\ \ \hbox{in}\ \ \Omega,\]
and let $U_{\lambda}$ be given by (\ref{ULambda}). Suppose $\phi \in C^1(\overline{\Omega})$ and $\phi \equiv C$ on $\partial \Omega$. 
Define the equimeasurable symmetric rearrangement $\phi^*$ of $\phi$, with respect to the measures $e^{w}dy$ and $e^{U_\lambda}dy$, by (\ref{rearrang}). Then $j^*(t), J^*(t)  $ are  absolutely continuous and decreasing  in $t>\min_{y\in \overline \Omega}\phi$  and
\[\int_{\{\phi^*=t\}}|\nabla \phi^*| ds \leq \int_{\{\phi=t\}}|\nabla \phi|ds,  \ \hbox{for  \it {a.e.}}  \ \ t>\min_{y\in \overline \Omega}\phi. \]
\end{proposition}

\section{The Sphere Covering Inequality}

The main objective of this section is to prove the following theorem. 

\begin{theorem}\label{MainTheorem}
Let $\Omega $ be a simply-connected  subset of $R^2$ and assume $w_i \in C^2(\overline{\Omega})$, $i=1,2$ satisfy
\begin{equation}\label{MainEq}
\Delta w_i +e^{w_i}=f_{i}(y),
\end{equation}
where  $f_2 \geq  f_1 \geq 0$ in $\Omega$. If $w_2 \ge w_1, w_2 \not \equiv w_1$ in $\Omega$ and $w_2=w_1$ on $\partial \Omega$, then
\begin{equation}
\int_{\Omega} e^{w_1}+e^{w_2}dy \geq 8\pi. 
\end{equation}
Moreover,  the equality only holds when $f_2 \equiv f_1 \equiv 0$ and there is an analytic function $\xi=g(y)$  such that 
$g(\Omega)=B_R$ and $g(\partial \Omega)= \partial B_R$  for some $R>0$,   and $e^{w_1(y)} dy= e^{U_{\lambda_1}(\xi)} d\xi,
e^{w_2(y)} dy= e^{U_{\lambda_2} (\xi)} d\xi $ with $ e^{U_{\lambda_1}(R)}=e^{U_{\lambda_2}(R)}=1$.
\end{theorem}

\begin{remark}
Note that  if $ f_1 \equiv f_2$,  the condition  $w_2 \ge w_1, w_2 \not \equiv w_1$  in the theorem  can just be replaced by $ w_2 \not \equiv w_1$, since there must be a simply connected subset $\omega$  of $\Omega$  such that  $w_2 >w_1$ in  $\omega$ and $w_2=w_1$ on $\partial \omega$ after  switching the indices.   If $w_2-w_1$ changes sign in $\Omega$,  the inequality has indeed a lower bound as $16 \pi$. 
\end{remark}

Before proving the above theorem, let us first show that Theorem \ref{MainTheorem} holds when $w_1 , w_2$ are both radial. Choose $\lambda_2>\lambda_1$ and let $U_{\lambda_1} , U_{\lambda_2}$ be given by (\ref{ULambda}).  Suppose $U_{\lambda_1}=U_{\lambda_2}$ on $\partial B_R$ for some $R>0$.  Then 
\[\frac{\lambda_1}{1+\frac{\lambda_1^2R^2}{8}} =\frac{\lambda_2}{1+\frac{\lambda_2^2R^2}{8}}=\kappa.\]
Hence $\lambda_1,\lambda_2$ are positive real roots of the quadratic equation 
\[R^2 \lambda^2+8=\frac{8}{\kappa} \lambda. \]
This implies $\kappa \leq 2/R^2$,   
\begin{equation}\label{scaling}
\lambda_1+\lambda_2=\frac{8}{\kappa R^2}, \ \ \hbox{and} \ \ \lambda_1 \lambda_2=\frac{8}{R^2}.
\end{equation}
Direct computations yield 
\begin{eqnarray*}
\int_{B_R}e^{U_{\lambda_1}}+e^{U_{\lambda_2}}dy &=& 8 \pi (\frac{\lambda_1^2R^2 }{8+\lambda_1^2R^2}+\frac{\lambda_2^2R^2 }{8+\lambda_2^2R^2})\\
&=& 8 \pi (\frac{\lambda_1^2R^2 }{\frac{8 \lambda_1}{\kappa}}+\frac{\lambda_2^2R^2 }{\frac{8 \lambda_2}{\kappa}})=8 \pi [\frac{\kappa R^2}{8}(\lambda_1+\lambda_2)]\\
&=& 8\pi. 
\end{eqnarray*}
Thus we have the following
 
\begin{proposition} \label{SymmetConjec}

Let $\lambda_2>\lambda_1$, and $U_{\lambda_1}$ and $U_{\lambda_2}$ be radial solutions of the equation
\[\Delta u +e^{u}=0, \]
defined in (\ref{ULambda}) with $U_{\lambda_2}>U_{\lambda_1}$ in $B_{R}$, and $U_{\lambda_1}=U_{\lambda_2}$ on $\partial B_{R}$, for some $R>0$. Then 
\[\int_{B_{R}}(e^{U_{\lambda_1}}+e^{U_{\lambda_2}})dy=8 \pi. \]
\end{proposition}

To understand the above inequality geometrically, we may scale the conformal factor by $1/2$ and consider two surfaces $\calS_1$ and $\calS_2$ with constant Gaussian curvature 1 as follows
\[ \calS_{1}=(B_R, e^{2V_{\lambda_1}}dy) \ \ \hbox{ and} \ \  \calS_{2}=(B_R, e^{2V_{\lambda_2}}dy).\]
where $2V_\lambda =U_\lambda-\ln2$.
Notice that the metrics $g_{i}=e^{2V_{\lambda_i}}dy$ have the same conformal factor on $\partial B_R$ and hence \eqref{scaling} holds and
\begin{equation}\label{roots}
\kappa(\frac{1}{\lambda_1}+\frac{1}{\lambda_2})=\frac{\kappa(\lambda_1+\lambda_2)}{\lambda_1 \lambda_2}=1.
\end{equation}

 Next we explain  that areas of $\calS_{1}$ and $\calS_{2}$ are equal to the areas of two complimenting spherical caps on the unit sphere, and consequently the total area must be 
\[A_{1}+A_{2}=\int_{B_R}e^{2V_{\lambda_1}}dy+\int_{B_R} e^{2V_{\lambda_2}}dy=4 \pi. \]
Indeed, we have 
\begin{eqnarray*}
\int_{B_R}e^{2V_{\lambda_i}}dy&=&\int_{B_R} \frac{\lambda_i^2}{2(1+\frac{\lambda_i^2|x|^2}{8})} dy\\
&=&\int_{B_{\frac{\lambda_i R}{\sqrt{8}}}}\frac{4}{(1+|x|^2)}dy, \ \ i=1,2. 
\end{eqnarray*}
Hence by the stereographic projection, $A_{i}$ is equal to the area of the spherical cap $\calC_i$ of the unit sphere  that lies below the plane 
\[z=z_i:=\frac{\frac{\lambda_i^2 R^2}{8}-1}{\frac{\lambda_i^2 R^2}{8}+1}=\frac{\frac{\lambda_i}{\kappa}-2}{\frac{\lambda_i}{\kappa}}=1-\frac{2\kappa}{\lambda_i}, \  \  i=1,2.\]
It follows from (\ref{roots}) that  $z_1<0$ while $z_2>0.$ Moreover  
\begin{eqnarray*}
|z_1|= \frac{2\kappa}{\lambda_1}-1=1-\frac{2\kappa}{\lambda_2}=z_2, \ \ i=1,2. 
\end{eqnarray*}
Thus $\calC_1$ and $\calC_2$ are two complimenting spherical caps on the unit sphere.  Note that the area of the smaller cap $\calC_1$ can be arbitrarily close to 0 or $2\pi$.

\vspace{.5cm}
The following lemma will play a key role in the proof of Theorem \ref{MainTheorem}. 

\begin{lemma}\label{LastEstimate}
Assume that $\psi \in C^{0, 1} (\overline{B_R})$ is an strictly decreasing, radial, Lipschitz function,  and  satisfies 
\begin{equation}\label{superSol}
\int_{\partial B_r} |\nabla \psi|  \le \int_{B_r}e^{\psi}
\end{equation}
{\it a.e.}  $r\in (0,R)$ and $\psi=U_{\lambda_1}=U_{\lambda_2}$  for some $\lambda_2> \lambda_1$ on $\partial B_R$, and $R>0$. Then 
there holds 
\begin{equation}\label{lastEstimate}
\hbox{either }  \,\,  \int_{B_R} e^{\psi} \le  \int_{B_R} e^{U_{\lambda_1}}  \,\,\quad  \hbox{or} \,\, \quad  \int_{B_R} e^{\psi}\geq   \int_{B_R} e^{U_{\lambda_2}}.
\end{equation}
Moreover if the inequality in (\ref{superSol}) is strict in a set with positive measure in $(0,R)$, then the inequalities in (\ref{lastEstimate}) are also strict.
\end{lemma}

{\bf Proof.} Let $m_1:=\int_{B_R}e^{U_{\lambda_1}}$, $m_2:=\int_{B_R}e^{U_{\lambda_2}}$, and $m:=\int_{B_R}e^{\psi}$. Also define 
\[\beta:=\left( \int_{\partial B_R} e^{\frac{\psi}{2}}\right)^2=\left( \int_{\partial B_R} e^{\frac{U_{\lambda_1}}{2}}\right)^2=\left( \int_{\partial B_R} e^{\frac{U_{\lambda_2}}{2}}\right)^2.\]
It follows from Proposition \ref{RadialBolProp} that 
\[\beta\geq \frac{1}{2}m(8 \pi -m).\]
On the other hand, 
\[\beta=\frac{1}{2}m_1(8 \pi -m_1)=\frac{1}{2}m_2(8 \pi -m_2).\]
Hence $m_1$ and $m_2$ are roots of the quadratic equation
\[x^2-8\pi x +2 \beta=0.\]
Since $m$ satisfies
\[m^2-8\pi m +2\beta\geq 0,\]
we have 
\[ \hbox{either} \quad m \le m_1  \quad \hbox{or } \quad m\ge m_2 .\]
Equality  holds only  when the equality in \eqref{superSol} holds for {\it a.e.} $r \in (0, R).$
This completes the proof.  \hfill $\Box$ \\ 

Now we are ready to prove Theorem \ref{MainTheorem}. \\ \\
{\bf Proof of Theorem \ref{MainTheorem}.} Suppose $w_1$ and $w_2$ satisfy the assumptions of Theorem \ref{MainTheorem}. Then 
\[\Delta (w_2-w_1)+e^{w_2}-e^{w_1}=f_2-f_1\geq 0.\]
Now we can choose $\lambda_2>\lambda_1$ such that  $U_{\lambda_1}$ and $U_{\lambda_2}$ are as described in Proposition \ref{SymmetConjec}, and 
\begin{equation}
\int_{\Omega}e^{w_1}=\int_{B_1}e^{U_{\lambda_1}}.
\end{equation}

Let $\varphi$ be the symmetrization of $w_2-w_1$ with respect to the measures $e^{w_1}dy$ and $e^{U_{\lambda_1}}dy$. 
Then by Proposition \ref{rearrangProp}
, 
\begin{eqnarray*}
\int_{\{\varphi=t\}}|\nabla \varphi| &\leq &  \int_{\{w_2-w_1=t\}}|\nabla (w_2-w_1)| \\
&\leq&\int_{\Omega_t} \bigl(e^{w_2}-e^{w_1} \bigr)\\
&=& \int_{\{\varphi>t\}}e^{U_{\lambda_1 }+\varphi}- \int_{\{\varphi>t\}}e^{U_{\lambda_1}}\\
&=& \int_{\{\varphi>t\}}e^{U_{\lambda_1}+\varphi}-\int_{\{\varphi=t\}} |\nabla U_{\lambda_1}|,  \quad  \hbox{for  {\it a.e.} }  t>0.
\end{eqnarray*} 
  Hence
\begin{equation}\label{Sub}
\int_{\{\varphi=t\}} |\nabla (U_{\lambda_1}+\varphi)| \le \int_{\{\varphi>t\} } e^{(U_{\lambda_1}+\varphi)},  \quad  \hbox{for  {\it a.e.} }  t>0.
\end{equation}
  Since  $\varphi \geq 0$ is decreasing in $r$,     $\psi:= U_{\lambda_1}+\varphi$ is a strictly decreasing function, and 
\begin{equation}\label{supersolution}
\int_{\partial B_r} |\nabla \psi| \le  \int_{B_r} e^{\psi} dy, \quad {\it a.e.}  \quad r \in (0, 1), 
\end{equation}
by Proposition \ref{rearrangProp} and the above inequality we know that  $\psi$ belongs to  $C^{0, 1} (B_R)$.  It follows from Lemma \ref{LastEstimate} that 
\[\int_{B_1}e^{\psi}dy =\int_{B_1}e^{U_{\lambda_1}+\varphi}dy \geq \int_{B_1}e^{U_{\lambda_2}}.\]
Hence
\[\int_{\Omega}e^{w_1}+e^{w_2}dy= \int_{B_1}(e^{U_{\lambda_1}}+e^{U_{\lambda_1}+\varphi})dy\geq \int_{B_1}e^{U_{\lambda_1}}+e^{U_{\lambda_2}}dy=8\pi.\]
Moreover,  it is clear that if  the equality holds,  then  $f_2 \equiv f_1$ and  the equality in \eqref{supersolution} holds for all $ r >0$.
This leads to the equality in Bol's inequality for $w_2$ in $\Omega$,  therefore $f_2 \equiv 0$  and there is an analytic function $\xi=g(y)$  such that  $g(\Omega)=B_R$ and $g(\partial \Omega)= \partial B_R$  for some $R>0$ and $e^{w_1(y)} dy= e^{U_{\lambda_1}(\xi)} d\xi$.  
Furthermore,  $\psi=U_{\lambda_1}+\varphi  \equiv U_{\lambda_2}(\xi)$.   This proof is complete.   \hfill $\Box$ 
 
Note that  the following consequence of  Bol's inequality and the equimeasurable symmetric rearrangement (see  Lemma 3.1 in \cite{GL-MR2670931} or  Proposition 3.3 in \cite{LL-MR2352519} for a proof) may be regarded as a limiting case of the Sphere Covering Inequality. 

\begin{proposition}\label{4piBound}
Let $\Omega \subset \R^2 $ be a simply-connected domain and  assume that $w\in C^2(\overline{\Omega})$ satisfies 
$\Delta w+e^{w} \ge 0$ in $\overline{\Omega}$ and $\int_{\Omega} e^w \le 8 \pi$.   Consider an open set $\omega \subset \Omega $  and define
the first eigenvalue of the operator $ \Delta + e^w $ in $H_0^1(\omega) $ by
$$ 
\lambda_{1, w} (\omega):=  \inf_{\phi \in H^1_0 (\omega) }  \bigl(\int_{\omega} |\nabla \phi| ^2 - \int_{\omega} \phi^2 e^w \bigr)   \le 0.
$$
Then $\int_{\Omega}e^{w} \ge 4 \pi$  if $\lambda_{1, w} (\Omega) \le 0$.
\end{proposition}

Suppose that  $w_1^k, w_2^k $  are solutions of  \eqref{MainEq}  with  $f_1^k,  f_2^k $ in  
$\Omega_k,   k=1, 2, \cdots$  and   the conditions of Theorem \eqref{MainTheorem}  hold for each $k$ large. Further assume that   $\Omega_k \to \Omega  \text{ in }  C^2,  w_i^k \to w  \text{ in } C^2,  f_i^k \to  f,   \text{ in }  C^0, i=1, 2$ as $ k \to \infty$.   Then  $w$ satisfies the condition in Proposition \ref{4piBound},  and the Sphere Covering Inequality  gives 
  the  same conclusion as Proposition \ref{4piBound}.  
  
  \begin{remark}\label{Remark}
  It can be seen from its  proof that Proposition \ref{4piBound} has a geometrical interpretation as follows:  Given a simply connected region $\omega$  on  a surface with Gaussian curvature less than $1/2$,   if the first eigenvalue  of Laplacian in $\omega$ with zero Dirichlet  boundary condition
   $$
   \lambda_1(\omega) :=\inf_{\phi \in H^1_0 (\omega) } \frac{ \int_{\omega} |\nabla \phi| ^2 }{\int_{\omega} \phi^2 e^w } \le 1,
   $$
then  the area of $\omega$  must be  bigger than or equal to $4\pi$,  which is the area of  a hemisphere  with Gaussian curvature $1/2$.  Note that such a hemisphere has first eigenvalue equal to $1$ as the height  function from the  boundary equator of the hemisphere is the first eigenfunction.  In other words,  Proposition \ref{4piBound} is an immediate consequence of the extremal  eigenvalue 
theorem  which says that a  geodesic disc on the sphere  achieves the smallest first eigenvalue of Laplacian among all  surfaces with the same area and  the same Gaussian curvature upper bound.   It would be very interesting to see if there is a simple  but deep geometric explanation of Theorem \ref{MainTheorem}.
\end{remark}

\section{Best Constant in a Moser-Trudinger Type Inequality} 

Let us consider the  functional   $J_\alpha (u)$  defined in  \eqref{MTAO} in 
 $$\mathcal{M}:=\{u\in H^1(S^2): \int_{S^2}e^ux_j=0 \ \ \hbox{for}\ \ j=1,2,3\}.$$
 
In this section we shall prove that $\inf_{u\in \mathcal{M}} J_{\alpha}(u)=0$ for $\alpha \geq \frac{1}{2}$.  Critical points of $J_\alpha(u)$, up to an additive constant, satisfy 
\begin{equation}\label{PDE-on-Sphere}
\frac{\alpha}{2}\Delta u+e^u-1=0 \ \ \ \ \hbox{on}\ \ S^2. 
\end{equation}

Following \cite{GL-MR2670931}, let $\Pi$ be the stereographic projection $S^2\rightarrow \R^2$ with respect to the north pole $N=(1,0,0)$: 

\[\Pi:= \left( \frac{x_1}{1-x_3}, \frac{x_2}{1-x_3}\right).\]
Suppose $u$ is a solution of $(\ref{PDE-on-Sphere})$, and let
\[\bar{u}(y):=u(\Pi^{-1}(y)) \ \ \hbox{for}\ \ y\in \R^2.\]
Then $\bar{u}$ satisfies 
\begin{equation}
\Delta \bar{u}+\frac{8}{\alpha(1+|y|^2)^2}(e^{\bar{u}}-1)=0 \ \ \hbox{in} \ \ \R^2. 
\end{equation}
Now if we let 
\begin{equation}\label{v}
v=\bar{u}-\frac{2}{\alpha} \ln (1+|y|^2)+\ln (\frac{8}{\alpha}),
\end{equation}
then $v$ satisfies 

\begin{equation}\label{PDE-on-Plane}
\Delta v+(1+|y|^2)^{2(\frac{1}{\alpha}-1)}e^{v}=0 \ \ \hbox{in} \ \ \R^2,
\end{equation}
and 
\begin{equation}\label{Total}
\int_{R^2}(1+|y|^2)^{2(\frac{1}{\alpha}-1)}e^{v} dy=\frac{8\pi}{\alpha}. 
\end{equation}

\subsection{Uniqueness of Axially Symmetric Solutions}

For the convenience of the reader, we first use a new method to prove Conjecture A  for axially symmetric functions, which was originally proven in \cite{GW-MR1760786} and \cite{Lin1-MR1770683}.

\begin{lemma}\label{axial}
Let $\alpha \geq \frac{1}{2}$ and $u \in \mathcal{M}$ be a solution of (\ref{PDE-on-Sphere}).  If $u$ is axially symmetric, then $u\equiv 0$. 
\end{lemma}

{\bf Proof.} We may assume that $u$ is symmetric about $x_3$-axis, i.e. $u=g(x_3)$, $x_3 \in [-1,1]$. Since $\int_{S^2}e^u x_3 d\omega=0$,  $g$ could not be monotone in $x_3$ unless it is identically equal to a constant $C$. Therefore, if $u \not \equiv C$, then it must take either its absolute minimum or absolute maximum at some point $x_3^0 \in (-1,1)$. Without loss of generality we can assume  $x_3^0 \geq 0$ and $g(x_3^0)= \max_{[-1, 1]} g(x_3)$.   Now choose some point $p=(0,x_2^p ,x_3^p) \in S^2$ with $x_3^0< x_3^p <1$ and let $u_p(x)=u(R^{-1}(x) )$ for some $R\in SO(3)$ with $R(p)=(0,0,1)$.  Define $\bar{u}_p=u_p(\Pi ^{-1})$ and let

\[
v_p=\bar{u}_p-\frac{2}{\alpha} \ln (1+|y|^2)+\ln (\frac{8}{\alpha}).
\]
the $v_p$ satisfies \eqref{PDE-on-Plane} and \eqref{Total}. 
Now let 
\[\varphi_p (y):=y_2 \frac{\partial v_p}{\partial y_1}-y_1 \frac{\partial v_p}{\partial y_2}.\]
Note that  the set of critical points of $\bar{u}_p$ contains a closed simple curve $\calC \subset \R^2$  which contains the origin in its interier,  and  $\bar{u}_p$ takes its absolute  maximum on $\calC$.  On the other hand $v_p$ and $\varphi_p$ are evenly symmetric about $y_1$-axis, therefore 
\[\calC \cup \{y=(y_1,0): y \in \R^2\} \subset \varphi^{-1}_p(0).\]
Hence $\varphi^{-1}_p(0)$ divides $\R^2$ into at least four simply-connected  regions $\Omega_i$, $i=1,2,3,4$. Now let $w_p:=\ln((1+|y|^2)^{2(\frac{1}{\alpha}-1)}e^{v_p} )$. Then $w_p$ satisfies 
\[\Delta w_p+e^{w_p}=\frac{8(\frac{1}{\alpha}-1)}{(1+|y|^2)^2}>0 \ \ \hbox{in} \ \ \R^2.\]
On the other hand $\varphi_p$ satisfies
\[\Delta \varphi_p+e^{w_p}\varphi_p=0 \ \ \hbox{in} \ \ \R^2.\]
Thus it follows from (\ref{4piBound}) that 
\begin{eqnarray*}
\frac{8\pi}{\alpha}=\int_{\R^2}(1+|y|^2)^{2(\frac{1}{\alpha}-1)}e^{v_p} dy=\sum_{i=1}^{4}\int_{\Omega_i}e^{w_p}>4\time 4 \pi=16 \pi.
\end{eqnarray*}
This implies $\alpha<\frac{1}{2}$ which is a contradiction. Therefore $\varphi_p \equiv 0$ and consequently $u$ is also axially symmetric about the line passing through $p$ and the origin. Since $p\neq (0,0,1)$, $u$ must be identically equal to a constant, and therefore must be zero.  $\hfill$ $\Box$\\

Next we prove that if $u$ is evenly symmetric about a plane passing through the origin, then $u$ is axially symmetric.  Note that this result was remarked by Ghoussoub and Lin \cite{GL-MR2670931},  we provide the details here since it is needed in the proof of the main result. 

\begin{lemma}\label{symmetry-wrp-plane}
Let $\alpha \geq \frac{1}{2}$ and $u$ be a solution of (\ref{PDE-on-Sphere}). If $u$ is evenly symmetric about a plane passing through the origin, then $u$ is axially symmetric.
\end{lemma}
{\bf Proof.}  The proof is similar to the proof of Lemma \ref{axial}. We may assume that $u$ is evenly symmetric about $x_1x_3$-plane. Let $u_0$ be the restriction of $u$ to 
$\{x \in S^2: \ \ x_2=0\}$ and assume that $p \in S^2$ is a maximum point of $u_0$. Since $u$ is symmetric about $x_1x_3$-plane, $p$ is also a critical point of $u$ on $S^2$. Without loss of generality we may assume $p=(0,0,-1)$.  We claim that 
\[\varphi(y_1, y_2)=y_2 \frac{\partial v}{\partial y_1}-y_1\frac{\partial v}{\partial y_2}\equiv 0,\]
where $v$ is defined by (\ref{v}). Suppose $\varphi \not \equiv 0$. Since $v$ has a critical point at the origin, the nodal line of $\varphi$ divides a neighborhood of the origin into at least four regions. On the other hand $\varphi$ is symmetric with respect to $y_1$-axis and the nodal line of $\varphi$ contains the $y_1$-axis. Therefore the nodal line of $\varphi$ divides $\R^2$ into at least 4 simply-connected  regions $\Omega_i$, $i=1,2,3,4$. 
As before, we can show that  $\varphi \equiv 0$ and consequently $u$ is axially symmetric about the line passing through $p$ and the origin. \hfill $\Box$

\subsection{The General Case} 

We  shall prove the even symmetry of a solution to \eqref{PDE-on-Sphere}. 

\begin{theorem}\label{GeneralCaseTheorem}
Let $\alpha \geq \frac{1}{2}$ and assume $u$ be a solution of (\ref{PDE-on-Sphere}). Then $u$ is  evenly symmetric about  any plane passing through the origin and a critical point of $u$. Therefore  $u$ must be  axially symmetric and consequently $u\equiv 0$. 
\end{theorem}

{\bf Proof.} Without loss of generality we may assume that $(1,0,0)$ is a critical point of $u$,   and  that $u$ is not symmetric about $x_1x_2$-plane. To finish the proof, it is enough to prove that $u$ is symmetric about the $x_1x_2$-plane. Define $u^*(x_1, x_2, x_3):=u(x_1, x_2,-x_3)$ and $\tilde{u}(x)=u(x)-u^*(x)$. Notice that $\tilde{u}(x_1, x_2, 0)=0$, for all $(x_1, x_2, 0)\in S^2$. Then $\tilde{u}$ satisfies 
\begin{eqnarray}
\frac{\alpha}{2}\Delta \tilde{u}+c(x)\tilde{u}=0, \ \ \hbox{on} \ \ S^2,
\end{eqnarray}
where 
\[c(x):=\frac{e^{u}-e^{u^*}}{u-u^*}.\]
Since $(1,0,0)$ is a critical point of $u$, it follows from the Hopf's lemma that $\tilde{u}$ must change sign in $S^+:=\{(x)\in S^2: x_3>0\}$. Therefore the nodal line of $\tilde{u}$ divides $S^{+}$ into at least two simply-connected  regions and there exists $S^{+}_+, S^{+}_- \subset S^+$ such that $u=u^*$ on $\partial (S^+_+ \cup S^+_-)$,  
\[u>u^* \ \ \hbox{on} \ \ S^+_+ \ \ \hbox{and} \ \ u<u^* \ \ \hbox{on} \ \ S^+_-.\]
Define $S^{-}_+, S^{-}_-$ to be the reflections of $S^{+}_+, S^{+}_-$ with respect to the $x_1x_2-$plane. Then we also have $u=u^*$ on $\partial (S^-_+ \cup S^-_-)$,  
\[u<u^* \ \ \hbox{on} \ \ S^-_+ \ \ \hbox{and} \ \ u>u^* \ \ \hbox{on} \ \ S^-_-.\]
Let $\Omega_1, \Omega_2,\Omega_3,\Omega_4 \subset \R^2$ be the images of $S^{-}_{-}, S^{-}_{+},S^{+}_{-},S^{+}_{+} \subset S^2$ under the stereographic projection, respectively. Define $v_1, v_2$ as follows
\[v_1(y)=u(\Pi^{-1}(y))-\frac{2}{\alpha}\ln(1+|y|^2)+\ln(\frac{8}{\alpha})\]
and
\[v_2(y)=u^*(\Pi^{-1}(y))-\frac{2}{\alpha}\ln(1+|y|^2)+\ln(\frac{8}{\alpha}).\]
Then $v_1$ and $v_2$ both satisfy (\ref{PDE-on-Plane}) and $w_i$ defined by 
\[w_i:=\ln((1+|y|^2)^{2(\frac{1}{\alpha}-1)}e^{v_i})\]
satisfies 
\[\Delta w_i+e^{w_i}=\frac{8(\frac{1}{\alpha}-1)}{(1+|y|^2)^2}\geq 0 \ \ \hbox{in} \ \ \R^2, \ \ i=1,2.\]
Moreover $w_1=w_2$ on $\partial \Omega_i$, $i=1,2,3,4$. Applying  the Sphere Covering Inequality (Theorem \ref{MainTheorem}) on $\partial \Omega_i$, $i=1,2,3,4$, 
we obtain that  

\begin{eqnarray*}
2 \times \frac{8\pi}{\alpha} &=&\int_{\R^2}(1+|y|^2)^{2(\frac{1}{\alpha}-1)}e^{v_1}dy+\int_{\R^2}(1+|y|^2)^{2(\frac{1}{\alpha}-1)}e^{v_2}dy\\
&\geq & \sum_{i=1}^{4}\int_{\Omega_i}e^{w_1}+e^{w_2}dy> 4\times 8 \pi.
\end{eqnarray*}
Hence $\alpha < \frac{1}{2}$ which is a contradiction. Thus $u$ is evenly symmetric about the $x_1x_2$-plane and the proof is complete. \hfill $\Box$

\section{Radial Symmetry of Solutions in $\R^2$}
In this section,  we shall consider solutions to a general class of equations in $\R^2 $  and prove radial symmetry  of the solutions. 
Assume  $u \in C^2 (\R^2) $  satisfies 
\begin{equation}\label{general}
\Delta u+k(|y|)e^{u}=0 \ \ \hbox{in}\ \ \R^2,
\end{equation}
and 
\begin{equation}\label{16Pi}
\frac{1}{2\pi} \int_{\R^2}k(|y|)e^u dy =\beta \leq 8,
\end{equation} 
where  $K(y)=k(|y|) \in C^2(\R^2)$ is  a non constant  positive function  satisfying 
\begin{eqnarray*}
&(K1) \quad \quad   &\Delta \ln(k(|y|))\geq 0,  \quad y \in \R^2\\
&(K2)  \quad \quad  &k(|y|) \le C(1+|y|) ^m,  \quad y \in \R^2
\end{eqnarray*}
for some constant $C, m>0$. It is easy to see that $(K1)$ implies  that  both  $k(r)$  and 
$\frac{r k'(r)}{k(r)} $ are  nondecreasing.  Let 
$$
2l=\lim_{r \to \infty } \frac{r k'(r)}{k(r)}.
$$
From  $(K2)$  we know that $ 0 \le 2l \le m$ and hence for any $\epsilon >0$ there exists a positive constant $C_\epsilon>0$ such that 
$$
C_\epsilon (1+|y|^2) ^{l-\epsilon} \le k(|y|) \le C(1+|y|^2) ^l,  \quad y\in \R^2. 
$$ 
Without loss of generality we may assume that $m=2l$.  
Then it follows from Theorem 1.1 in \cite{ChengLin-MR1446203} that 
$$\beta \ge 2l+2.$$ 
Following \cite{GL-MR2670931} and using Pohazaev identity,  we can obtain the following result. 
  
\begin{proposition}\label{range}
Suppose $u$ is a solution to \eqref{general}-\eqref{16Pi},  where  $(K1)-(K2)$  hold with $m=2l$.  Then,  if $ \beta > 2l+2$,  there holds 
\begin{equation}\label{range}
4 < \beta < 4l+4.
\end{equation}
\end{proposition}

{\bf Proof.} By Theorem 1.1 in \cite{ChengLin-MR1446203},  we have 

\begin{equation}\label{asymp}
u(y)=-\beta \ln (|y|)+ C+ O(|y|^{-\gamma}) 
\end{equation}
for some constants $C$ and $\gamma>0$ as $x \to \infty$, if $ \beta > 2l+2$. Also if $\beta =2l+2$, then for any $\epsilon>0$  there exists $ R(\epsilon) >0$ such that 
$$
-\beta \ln (|y|) -C \le u(y)\le  (\epsilon-\beta ) \ln (|y|),  \quad |y| \ge R(\epsilon)
$$
for some constant $C$. On the other hand, it is easy to see that when $\beta > 2l+2$,  we have
$$
\nabla u=( -\beta +o(1) )  \frac{x}{|y|^2} , \quad \text{ as }  y \to \infty.
$$

Multiplying \eqref{general} by $y \cdot \nabla u $ and integrating by parts on $B_R=\{ y:  |y| \le R\}$,  we obtain

\begin{eqnarray*}
&\int_{\partial B_R} (y\cdot \nabla u) \frac{\partial u}{\partial \nu} ds-\frac{1}{2} \int_{\partial B_R} (y\cdot \nu) |\nabla u|^2 ds=-\int_{ B_R} k(|y|) y \cdot \nabla e^u dy \\
&=\int_{ B_R} ( 2k(|y|) +k'(|y|) |y| ) e^u dy-\int_{\partial B_R} (y\cdot \nu) k(|y|) e^uds.
\end{eqnarray*} 
Letting $R \to \infty$ and using \eqref{asymp},  we obtain that 
$$
\int_{ \R^2} \bigl( 2k(|y|) +k'(|y|) |y| \bigr) e^u dy= \pi \beta^2.
$$
Hence we derive  \eqref{range} from 
$$   2 k(|y|) \le 2k(|y|) +k'(|y|) |y| \le (2l+2 ) k(|y|),  \quad y \in \R^2, 
$$
and the fact that  equality  holds in the above inequalities only when $l=0$  and $k$ equals to a constant.  Note that  by our assumptions, 
$k(|y|)=|y|^{2l}$ is not allowed for $l >0$ since $k(0)=0$ in this case, nor is  $k$  allowed to be equal to  a positive constant.  The proof is complete.  
\hfill $\Box$

\begin{remark}
 In all applications considered in this paper,  it holds that  $\beta >2l+2$.   We wonder if   $\beta > 2l+2$  is always true  for all solutions to \eqref{general}-\eqref{16Pi}  under the general conditions $(K1)-(K2)$.
\end{remark}

It is shown in \cite{Lin1-MR1770683} that 

\begin{proposition}\label{radial}
If  $0<l \le 1, $ then  there exists a  radially symmetric
solution $u_\beta$ to \eqref{general} if and only if   $\beta \in (4, 4l+4)$.  The radial solution is also unique in this case.  Also If $l>1$ and $\beta \in (4l, 4l+4)$, then there exists a unique radially symmetric solution $u_\beta$  to \eqref{general}.

\end{proposition}

Now we are ready to prove the following general theorem.

\begin{theorem} \label{th-general}
Assume that $K(y)=k(|y|)>0$  satisfies  $(K1)-(K2),$  and  $u$ is a solution to \eqref{general}-\eqref{16Pi}  with $ 2l+2< \beta \le  8$.  
Then $u$ must be  radially symmetric. 
\end{theorem}

{\bf Proof.} Since $\lim_{|y|\rightarrow \infty} u(y)=-\infty$, $u$ has a maximum point $p\in \R^2$. We first prove that $u$ is evenly symmetric about the line passing through the origin and $p$. In particular if $p=(0,0)$, then the following argument guarantees that $u$ is evenly symmetric about any line passing through the origin and hence $u$ must be  radially symmetric. Without loss of generality we may assume that $p$ lies on $y_1$-axis. Define 

\begin{equation}
v(y_1,y_2)=u(y_1,y_2)-u(y_1,-y_2). 
\end{equation}
Suppose $v\not \equiv 0$. Then the nodal line of $v$, $v^{-1}(0)$, contains the $y_1$-axis. On the other hand since the critical point $p$ lies on $y_1$-axis, the nodal line of $v$ divides every small neighborhood of $p$ into at least four regions. Therefore the nodal line of $v$ divides $\R^2$ into at least four simply-connected  regions $\Omega_i$, $i=1,2,3,4$.  Now notice that on each $\Omega_i$ the equation 
\[\Delta u+k(|y|)e^u=0 \ \ y \in \Omega_i\]
has two solutions $u_i^1(y_1,y_2)=u(y_1,y_2)$ and $u_i^2(y_1,y_2)=u(y_1,-y_2)$ with $u_i^1|_{\partial \Omega}=u_i^2|_{\partial \Omega_i}$. Define $w:=u+\ln(k(|y|))$. Then $w$ satisfies 
\begin{equation}\label{generalPrime}
\Delta w+e^{w}=\Delta(\ln(k(|y|)))\geq 0.
\end{equation}
Thus on each $\Omega_i$, the above equation has two solutions $w_i^1,w^2_i$ with $w^1_i|_{\partial \Omega_i}=w_i^2|_{\partial \Omega}$, $i=1,2,3,4$. Hence it follows from Theorem \ref{MainTheorem} that 
\begin{eqnarray*}
4\pi \beta&=& 2\int_{\R^2}k(|y|)e^udy=\int_{\R^2}k(|y|)e^{u(y_1,y_2)}dy+\int_{\R^2}k(|y|)e^{u(y_1,-y_2)}dy \\
&\geq & \sum_{i=1}^4 \int_{\Omega_i}e^{w_1}+e^{w_2}dy>4 \times 8\pi=32\pi . 
\end{eqnarray*}
Consequently $\beta> 8$ which is a contradiction, and therefore $u$ is evenly symmetric about the $y_1$-axis. 

Next we shall prove that $u$ is indeed axially symmetric. Let $\phi= y_2 \cdot u_{y_1}-y_1\cdot  u_{y_2} $.  Then $\phi$ satisfies
    \begin{equation}\label{phi}
      \Delta \phi + K(y) e^u \phi= 0, \quad y \in \R^2.
    \end{equation}
On the other hand,   $u_{y_2}$  satisfies 
     \begin{equation}\label{u2}
      \Delta u_{y_2} + K(y) e^u u_{y_2} = -y_2 \frac{k'(|y|) }{|y|} e^u, \quad y \in \R^2.
    \end{equation}
Note that both $u_{y_2}$ and $\phi$ are odd function in $y_2$. Let us multiply  equation (\ref{phi}) by $u_{y_2}$ and equation (\ref{u2}) by $\phi$ and subtract. Then,  integrating the resulting equation in $B_R^+=\{ y:  y_2>0,  |y| \le R\}$,  we obtain
    $$
    \int_{\partial B_R^+} \phi  \frac{\partial u_{y_2} }{\partial \nu} -u_{y_2} \frac{\partial \phi}{\partial \nu} ds = -\int_{B_R^+} y_2 \frac{k'(|y|) }{|y|} e^u \phi.
    $$
Applying the standard Schauder estimates for the elliptic equation  satisfied by $u+\beta \ln(|y|)$ and using the fact  that $ \beta > 2l+2$,  we obtain
    $$
    |\nabla u(y) | \le \frac{ C}{|y|},  \quad  \quad  |\nabla^2 u(y)|  \le  \frac{ C}{|y|^2}, \quad |y|>1
    $$
    for some constant $C$. Letting $R \to \infty$,  we  derive
    $$
    \int_{{\R^2}_+}  y_2 \frac{k'(|y|) }{|y|} e^u \phi dy =0.
    $$
    
   We claim that  $\phi \equiv 0 $ in $\R^2$.    Assume the contrary.   Since  $y_2 \frac{k'(|y|) }{|y|} e^u>0 $ in $ {\R^2}_+$,   there exist  at least  two regions $\Omega_1, \Omega_2 \subset {\R^2}_+$ such that $\phi> 0 $ in $\Omega_1$ and $\phi < 0 $ in $\Omega_2$  and $\phi=0 $ on   
   $\partial \Omega_i,  i=1, 2.$   Applying Proposition \ref{4piBound} to $\Omega_i,  i=1, 2$,    we  conclude that 
    $$
    \int_{{\R^2}^+}  k(|y|)  e^u  dy \ge \int_{\Omega_1}e^{w}dy+\int_{\Omega_1}e^{w}dy >8 \pi,
    $$
   and therefore $\beta >  8$.  This contradiction  shows $\phi \equiv 0$ in $\R^2$, and hence $u(y)$ is radially symmetric. \hfill $\Box$ \\ \\

Now we consider several special cases of the equation \eqref{general}-\eqref{16Pi}.  First,  if $K(y)= (1+|y|^2) ^l$ for some $ l \ge 0$,   then the equations \eqref{general}-\eqref{16Pi} read as 
\begin{equation}\label{l}
\Delta v+(1+|y|^2)^{l}e^{v}=0 \ \ \hbox{in} \ \ \R^2,
\end{equation}
and 
\begin{equation}\label{integral}
\frac{1}{2\pi} \int_{\R^2}(1+|y|^2)^{l}e^{v} dy=\beta. 
\end{equation}
The following result  is conjectured  in \cite{GL-MR2670931}. \\ \\
  {\bf Conjecture D.}   
For $ 0 < l \le 2$ and   $\beta=(2l+4)$,  solutions to \eqref{l}-\eqref{integral} must be radially symmetric.
\\ \\
For $ -2<l\le 0$,  the radial symmetry of solutions to \eqref{l}-\eqref{integral}  was shown in   \cite{CL1-MR1121147} and \cite{CK-MR1262195} by the moving plane method; while for  $l>0$ the moving plane method does not seem to work,  the conjecture was shown  
in \cite{GL-MR2670931} for $ 0 < l \le 1$ by using the Alexandrov-Bol inequality.  For $  2< l \not =(k-1)(k+2)$,  where $k\ge 2$ is an integer,  it is pointed out by Lin  in  \cite{Lin1-MR1770683} that  there is a non-radial solution to  \eqref{l}-\eqref{integral} .  A direct application of Theorem \ref{th-general}  to \eqref{l}-\eqref{integral}  leads to an affirmative answer to {Conjecture D}.  Indeed,  all solutions to \eqref{l}-\eqref{integral} must be radially symmetric as long as $ \beta \le 8$. 

Another example is the following equation from the study of self- gravitating strings for a massive W-boson model coupled to Einstein theory in account of gravitational effects  (\cite{PT-MR2876669}, \cite{Yang-MR1277473}).

 \begin{equation}\label{string}
\Delta v+(1+|y|^{2l}) e^{v}=0 \ \ \hbox{in} \ \ \R^2,
\end{equation}
and 
\begin{equation}\label{str-int}
\frac{1}{2\pi} \int_{\R^2}(1+|y|^{2l})e^{v} dy= \beta, 
\end{equation}
where $l>0$. It is shown in \cite{PT-MR2876669} that \eqref{string}-\eqref{str-int}  admit a radial solution if and only if 
\[\beta \in \bigl( 4 \max\{ 1, l \},  4(l+1) \bigr)\]
and  the corresponding radial solution is unique.  Furthermore, for $0 <l \le 1$,  the interval above is also optimal for the solvability of  \eqref{string}-\eqref{str-int} among  non-radial functions.  The main known difference between   \eqref{string} and \eqref{l} is that the latter possesses radial solutions for a larger range of $\beta$
 which is at least  $ (2l+2,  4l+4)$,  and has multiple radial solutions when $l> 2$ and $  \beta 
 \in (\beta_l, 4l)$ for some $\beta_l \in (2l+2, 2l+4)$, which also implies the existence of  non-radial solutions for $\eqref{l}$ for $l>2$  (see\cite{DEJ-MR3377875}, \cite{Lin1-MR1770683}).   While the former has  a radial solution only  for  $\beta \in \bigl( 4 \max\{ 1, l \},  4(l+1) \bigr)$,  which is  also unique.  In particular,   no non-radial solution is known in this case.
 
 Theorem \ref{th-general} implies that  solutions to \eqref{string}-\eqref{str-int}  must be radially symmetric when $ \beta \le 8$.   As a consequence,
 the solvability  range of $\beta$ among non-radial functions  must be $\beta >  4 \max\{ 1, l \}$ when $l \le 2$. \\

{\bf Proof of Theorem \ref{MeanFieldUniquenessTheorem}.} Let $u$ be a solution of (\ref{meanFieldPDE})-(\ref{meanFieldPDE1}). Without loss of generality we may assume $P=(0,0,1)$. Now let $\Pi: S^2 \rightarrow R^2$ be the stereographic projection with north pole at $P=(0,0,1)$.  Define 
\begin{eqnarray}
v(y)=u(\Pi^{-1} y)-\ln \left( \int_{S^2}e^{u} d\omega \right)+\ln(16\pi (3+\alpha))-3\ln(1+|y|^2),
\end{eqnarray}
where $y=\Pi(x)$. Then $u$ is a solution of (\ref{meanFieldPDE})-(\ref{meanFieldPDE1}) if and only if $v$ satisfies 
\begin{eqnarray*}
\left\{ \begin{array}{ll}
\Delta v+(1+|y|^2)e^{v}=0 &\text{in } \R^2\\
\int_{\R^2}(1+|y|^2)e^{v}dy=4\pi (\alpha+3).
\end{array} \right.
\end{eqnarray*}
Therefore it follows from Theorem \ref{th-general} that $v$ is radially symmetric about the origin. Hence $u$ is axially symmetric with respect to $\overrightarrow{OP}$ and the proof is complete.  \hfill $\Box$\\ \\

{\bf Proof of Theorem \ref{OnsagerVortexTheorem}.} Without loss of generality we may assume that $\vec{n} =(0,0,1)$. Let $\Pi: S^2 \rightarrow R^2$ be the stereographic projection with north pole at $\vec{n}=(0,0,1)$. Define 
\[v(y)=u(\Pi^{-1}(y)) \ \ \hbox{for} \ \ y\in \R^2.\]
Then $v$ satisfies 
\begin{equation}
\Delta v+\frac{J^2(y)\exp(\alpha v-\gamma \psi (y))}{\int_{\R^2}J^2(y)\exp(\alpha v-\gamma \psi (y)) dy }-\frac{J^2(y)}{4\pi}=0 \ \ \hbox{for}\ \ y\in \R^2, 
\end{equation}
where 
\[J(y)=\frac{2}{1+|y|^2} \ \ \hbox{and} \ \ \psi(y)=\frac{|y|^2-1}{|y^2|+1}.\]
Now define
\[w(y):=\frac{\alpha}{2} \left(v(y)-\frac{1}{4 \pi} \ln(1+|y|^2) \right)+c,\]
with
\[c=\frac{1}{2}\left( \gamma +\ln \left( \frac{2}{\alpha}\int_{\R^2} J^2(y)e^{\alpha v-\gamma \psi}\right) \right).\]
Then we have
\begin{equation}
\Delta w(y)+ K(y)e^{w}=0 \ \ \hbox{in} \ \ \R^2,
\end{equation}
and 
\[\int_{\R^2}K(y)e^{w(y)}dy=\alpha,\]
where 
\begin{equation}\label{onsager}
K(y)=2(1+|y|^2)^{(-2+\frac{\alpha}{4\pi})}e^{\gamma J(y)}.
\end{equation}
Now we compute 
\[\Delta(\ln{K(y)}) =\frac{4(-2+\frac{\alpha}{4 \pi})}{(1+|y|^2)^3}+\frac{8 \gamma(|y|^2-1)}{(1+|y|^2)^3}.\]
Since the right hand side of the above equation is nonnegative  for $0 \leq \gamma \leq \frac{\alpha}{8\pi}-1$, it follows from Theorem \ref{th-general} that $w$ is radially symmetric about the origin. \hfill $\Box$ \\

$\mathbf{Acknowledgement}$  The authors would like to thank Professors  Alice Chang,  Nassif Ghoussoub, Yanyan Li,  Fernando Marques, Richard Schoen,  Paul Yang for their interests and helpful comments on the earlier draft of the paper.   The first author is partially supported by a Simons Foundation Collaborative Grant (Award \#199305) and NSFC grant No 11371128.

\bibliographystyle{plain}
\bibliography{SphereCovering}

\end{document}